\documentclass[12pt]{amsart}
\usepackage{fullpage}
\usepackage{mathrsfs,eucal,amsmath,amsthm,amsfonts,amssymb,amscd,amsbsy,latexsym}
\usepackage[all,2cell]{xy} \UseAllTwocells \SilentMatrices
\usepackage{url}

\usepackage{verbatim}

\begin{document}


\newtheorem{teo}{}
\newtheorem{thm}{\textsc{Theorem}}

\newtheorem{dfn}[thm]{\textsc{Definition}}

\newtheorem{lem}[thm]{\textsc{Lemma}}

\newtheorem{cor}[thm]{\textsc{Corollary}}

\newtheorem{ex}[thm]{\textsc{Example}}

\newtheorem{prp}[thm]{\textsc{Proposition}}

\newtheorem{notation}[thm]{\textsc{Notation}}


\newcommand{\aroincl}{\xymatrix@!=0.01pc{\ar@{^{(}->}[r]&}}

\newcommand{\binomial}[2]{\left(\begin{array}{c} #1 \\ #2 \end{array}\right)}

\newcommand{\matriz}[4]{\left(\begin{array}{cc}#1 & #2 \\ #3 & #4\end{array}\right)}

\newcommand{\aro}{\xymatrix@!=0.01pc{\ar[r]&}}

\newcommand{\quadrado}[8]{\xymatrix{#1 \ar[r]^{#5}\ar[d]_{#8}&  #2 \ar[d]^{#6} \\   #4\ar[r]_{#7} &  #3}}

\newcommand{\quadradod}[4]{\xymatrix{#1 &  #2 \ar[l] \\   #4\ar[u] &  #3\ar[l]\ar[u]}}

\newcommand{\cart}{\ar@{}[dr]|{\square}}

\newcommand{\cartt}{\ar@{}[drr]|{\square}}

\newcommand{\rela}[4]{\xymatrix{#1  \ar@<+.7ex>[r]^{#2}\ar@<-.7ex>[r]_{#3} & #4}}

\newcommand{\seta}[3]{\xymatrix{#1\ar[r]^-{#2} & #3}}

\newcommand{\setaa}[3]{\xymatrix{#1\ar[rr]^{#2}&& #3}}

\newcommand{\arou}[1]{\xymatrix@!=0.01pc{\ar[r]^{#1}&}}




\newcommand{\mm}[1]{\mathrm{#1}}


\newcommand{\ccc}[1]{\mathcal{#1}}
\newcommand{\cc}[1]{\mathscr{#1}}

\newcommand{\g}[1]{\mathfrak{#1}}


\newcommand{\GA}{\mathbb{G}_a}

\newcommand{\GL}{\mathbb{GL}}

\newcommand{\NN}{\mathbb{N}}

\newcommand{\ZZ}{\mathbb{Z}}

\newcommand{\CC}{\mathbb{C}}

\newcommand{\QQ}{\mathbb{Q}}

\newcommand{\RR}{\mathbb{R}}

\newcommand{\FF}{\mathbb{F}}

\newcommand{\DD}{\mathbb{D}}

\newcommand{\qq}{\overline{\QQ}}

\newcommand{\PP}{\mathbb{P}}

\newcommand{\af}{\mathbb{A}}

\newcommand{\al}{\alpha}

\newcommand{\be}{\beta}

\newcommand{\ga}{\gamma}
\newcommand{\Ga}{\Gamma}

\newcommand{\om}{\omega}
\newcommand{\Om}{\Omega}

\newcommand{\te}{\theta}
\newcommand{\Te}{\Theta}

\newcommand{\ph}{\varphi}
\newcommand{\Ph}{\Phi}

\newcommand{\ps}{\psi}
\newcommand{\Ps}{\Psi}

\newcommand{\ep}{\varepsilon}

\newcommand{\vr}{\varrho}

\newcommand{\de}{\delta}
\newcommand{\De}{\Delta}

\newcommand{\la}{\lambda}
\newcommand{\La}{\Lambda}

\newcommand{\ka}{\kappa}

\newcommand{\si}{\sigma}
\newcommand{\Si}{\Sigma}

\newcommand{\ze}{\zeta}



\newcommand{\ot}{\otimes}
\newcommand{\ov}[1]{\overline{#1}}

\newcommand{\fr}[2]{\frac{#1}{#2}}

\newcommand{\lid}{\varinjlim}

\newcommand{\lip}{\varprojlim}
\newcommand{\hh}[3]{\mm{Hom}_{#1}(#2,#3)}

\newcommand{\vs}{\vspace{0.3cm}}

\newcommand{\pr}{\textbf{Proof:}}

\newcommand{\pd}{\partial}

\newcommand{\un}[1]{\underline{#1}}

\newcommand{\hs}{\hspace{1.5cm}}

\newcommand{\spc}[1]{\mm{Spec}\,#1}

\newcommand{\tyy}[3]{\left\{\begin{array}{c}#1 \\ #2 \\ #3 \end{array}\right}

\newcommand{\ty}[2]{\left\{\begin{array}{c}#1 \\ #2  \end{array}\right.}


\newcommand{\rep}[2]{\mm{Rep}_{#1}(#2)}

\newcommand{\strl}{\mathbf{str}^\#(\g{X/o})}
\newcommand{\strk}{\mathbf{str}(\g{X}/K)}
\newcommand{\stro}{\mathbf{str}(X_0/k)}
\newcommand{\str}{\mathbf{str}(\g{X/o})}
\newcommand{\go}{\g{o}}
\newcommand{\cm}{\cc{M}}
\newcommand{\repl}[2]{\mm{Rep}^{\#}_{#1}(#2)}
\newcommand{\bb}[1]{\mathbf{#1}}
\newcommand{\modulef}[1]{\text{$#1$-$\mm{mod}$}}
\newcommand{\module}[1]{\text{$#1$-$\mm{Mod}$}}
\newcommand{\lidd}[1]{\lid_{#1}}
\newcommand{\lipp}[1]{\lip_{#1}}

\title{On the rank of the fibers of elliptic K3 surfaces}
\author{Cecilia Salgado}
\address{Mathematisch Instituut Leiden, Niels Bohrweg 1}
\email{salgado@math.leidenuniv.nl}
\maketitle

\begin{abstract}
Let $X$ be an elliptic K3 surface endowed with two distinct Jacobian elliptic fibrations $\pi_i$, $i=1,2$, defined over a number field $k$. We prove that there is an elliptic curve $C\subset X$ such that the generic rank over $k$ of $X$ after a base extension by $C$ is strictly larger than the generic rank of $X$. Moreover, if the generic rank of $\pi_j$ is positive then there are infinitely many fibers of $\pi_i$ ($j\neq i$) with rank at least the generic rank of $\pi_i$ plus one. 
\end{abstract}
\maketitle

\section{Introduction}
 We investigate elliptic K3 surfaces $X$ over a number field $k$. By this, we mean a smooth projective geometrically connected surface $X$ with trivial canonical class and $H^1(X,\mathcal{O}_X)=0$ (\emph{K3 surface}) and with a flat morphism $\pi$ to a projective curve $B$ such that all but finitely many fibers are elliptic curves (\emph{elliptic surface}). Moreover, the elliptic fibrations considered throughout this paper are \emph{Jacobian}, i.e., endowed with a section defined over $k$.
  
 
  We are interested in comparing the Mordell-Weil rank $r_t$ of the fibers
 $\pi^{-1}(t)$ above a $k$-rational point $t\in B(k)$ with
 the Mordell-Weil rank $r$ of the generic fiber $X_{\eta}$.

  A theorem on specializations by N\'eron \cite{N}, treated in the setting of elliptic surfaces by Silverman \cite{Sil}, assures that for all but
 finitely many $t \in B(k)$, we have \[r_t \geq
 r.\]
 A natural question arises:

 \textbf{Question 1:} What can be said about the set
 \begin{equation}\label{conj}\{t \in B(k),
 r_t \geq r+1 \}?
\end{equation}

   Different methods have been used to tackle this problem for special classes of elliptic surfaces, but, so far, a general technique for solving the question in complete generality is lacking.
     
   If the fibration $X \rightarrow B$ is non-isotrivial, then Helfgott, studying the behavior of root numbers in families of elliptic curves, showed, under classical conjectures, that the above set is infinite (see \cite{He} and \cite{He2}). Rohrlich in \cite{R} and Manduchi in \cite{Ma} also analyzed root numbers in families of elliptic curves, obtaining results for elliptic surfaces in the case $k \simeq \QQ$ and $B \simeq \mathbb{P}^1$.
  
    Question 1 has not been explicitly treated so far in the literature for surfaces with Kodaira dimension zero or larger. In \cite{BogoTschinkel}, Bogomolov and Tschinkel proved potential density of the set of rational points on K3 surfaces endowed with an elliptic fibration or with an infinite automorphism group. As a corollary of their result, K3 surfaces with a rank zero elliptic fibration defined over a number field $k$ have infinitely many fibers with positive rank over a finite extension of $k$.

   In this paper, we prove that the set in (\ref{conj}) is infinite for K3 surfaces endowed with two distinct Jacobian elliptic fibrations with positive rank. This will follow from Theorem \ref{t1}, which states that given such a K3 surface with Jacobian elliptic fibrations $\pi_1$ and $\pi_2$, then the Mordell-Weil rank of $\pi_i$ jumps after a base change by an elliptic curve. This elliptic curve belongs to the set $\{\pi_j^{-1}(t); t \in \mathbb{P}^1 \}$, with $i\neq j$.

 \textbf{Motivation:} The Kodaira dimension of an elliptic surface can be -1, 0 or 1. Elliptic surfaces with Kodaira dimension -1 are rational elliptic surfaces or trivial ones, i.e., of the form $E \times \mathbb{P}^1$ where $E$ is an elliptic curve. For the latter, the rank of the fibers is constant, equal to the rank of $E$. The former class was the subject of \cite{B}, where Billard showed that if $k=\mathbb{Q}$ and $X$ is a non-isotrivial rational elliptic surface that is $\mathbb{Q}$-birational to $\mathbb{P}^2$, then the above set is infinite. Later the author extended the class of rational elliptic surfaces considering the problem over an arbitrary number field and showing that if $X$ is a rational elliptic surface with a $k$-minimal model of degree at least three, then the above set is infinite and, under some extra conditions, that the set 
\[ \{t \in B(k), r_t \geq r+2 \} \]  

 is also infinite (see \cite{eu} and \cite{salgado}).

   The class of surfaces with Kodaira dimension zero possibly admitting an elliptic fibration contains Enriques, K3 and bielliptic surfaces. For the first class, the fibration does not have a section since there are multiple fibers. Bielliptic surfaces are naturally equipped with two distinct elliptic fibrations. One of these fibrations has multiple fibers and thus no section. The other has an elliptic base $E$ and fibers generically isomorphic to $E$, in which case the problem of comparing the rank of the fibers becomes trivial. We are left with elliptic K3 surfaces. These surfaces cannot only admit an elliptic fibration with a section, but, moreover, they are the only surfaces that can be (Jacobian) elliptic in more than one way (admitting two distinct elliptic fibrations with a section).
   
 As a natural follow-up of \cite{salgado}, we address the problem of comparing the rank of the generic and special fibers for elliptic K3 surfaces obtaining a class of elliptic K3 surfaces for which the set in (\ref{conj}) is infinite (Theorem \ref{t1}). 
 
   Surfaces with Kodaira dimension 1 are always elliptic. Much less is known about their geometry and arithmetic.

\textbf{Remark 1:} The question raised above makes sense for any elliptic surface with a
 section. Certainly, due to Faltings' Theorem, the most interesting case concerns surfaces
 fibering over rational or elliptic curves $B$ with a positive Mordell-Weil rank.

 The paper is organized as follows. In Section 2, we give a very brief overview of elliptic K3 surfaces and state the main result of the paper (Theorem \ref{t1}). Section 3 is dedicated to proving this theorem, while Section 4 is devoted to examples of elliptic K3 surfaces with more than one elliptic fibration and to applying Theorem \ref{t1} to them.
 
\section{Elliptic K3 surfaces}

\subsection{Basic facts}
\indent\par

Let $X$ be an elliptic K3 surface. If the fibration $\pi: X \rightarrow B$ is Jacobian then $\mathrm{Pic}^0(X) \simeq \mathrm{Pic}^0(B)$ (notice that this holds for any elliptic surface). Thus $B \simeq \mathbb{P}^1$. 

Different from rational elliptic surfaces, whose Picard group has rank 10, K3 (not necessarily elliptic) surfaces satisfy:

    \[ \mathrm{rank}(\mathrm{Pic}(X)) \leq 20.\]
   
  In any elliptic surface the fiber components and the zero section generate a rank two subgroup of the Picard group and hence necessary condition for a K3 surface to have an elliptic fibration is $\mathrm{rank}(\mathrm{Pic}(X)) \geq 2$. For sufficiency, note that by \cite[Theorem 1]{Sha} $X$ has an elliptic fibration if and only if the Picard lattice $S_X$ represents zero, i.e., there exists $0 \neq x \in S_X$ such that $x^2=0$. This last condition is assured by a higher bound, namely \[\mathrm{rank}(\mathrm{Pic}(X)) \geq 5.\]

    Explicit examples of elliptic K3 surfaces for each possible value have been computed (see \cite{Kuwata} for elliptic K3 with Mordell-Weil rank between 0 and 18 except for 15. This last case was covered in \cite{Kloosterman}).     
         
     Another important characteristic of K3 surfaces is the possibility of admitting more than one Jacobian elliptic fibration. We will exploit this feature to give an answer Question 1.

 We now consider the subclass of elliptic K3 surfaces of our interest.

 \subsection{Two distinct elliptic fibrations with a section}
\indent\par
  We propose in this subsection the following question:

  \textbf{Question 2:} When is a K3 surface endowed with two distinct elliptic fibrations? 

\textbf{Example A:}  In the case $X$ is a Kummer surface associated to the product of two elliptic curves, Nishiyama described in \cite{Nishiyama} all possible distinct Jacobian elliptic fibrations in it. Kummer surfaces provide thus a good source of K3 surfaces with distinct elliptic fibrations on it.

\textbf{Example B:} An interesting class of elliptic K3 surfaces is the one of surfaces obtained by a quadratic base-change of a rational elliptic surface. Let $\pi_0: X_0 \rightarrow B_0$ be a rational elliptic surface and $\varphi: B\rightarrow B_0$ a degree two morphism. If $\varphi$ is not ramified above a non-reduced fiber of $\pi_0$ then $X= X_0 \times_{B_0}B$ is a K3 surface.
Under certain conditions these surfaces are also endowed with at least two distinct elliptic fibrations. Namely, if there is a non-constant linear pencil of rational curves $\{C_t\}_{t \in \mathbb{P}^1}$ in $X_0$ that are not in the same numerical class as $B$ and such that the morphism $C_t \rightarrow B_0$ given by the restriction of the fibration $X_0 \rightarrow B_0$ has degree two. The curves $C_t \times_{B_0}B$ give an elliptic fibration on $X$ that is clearly distinct from the fibration $X \rightarrow B$. 
Since $\mathrm{rank} (\mathrm{Pic}(X_0))=10$, $X$ will satisfy $\mathrm{rank}(\mathrm{Pic}(X))\geq 10$.

\textbf{Remark 2:}
To increase the rank of the fibers the new base must have infinitely many $k$-rational points. Thus, by Faltings' Theorem, we must have $g(C_t \times_{B_0} B)\leq 1$. An application of the Hurwitz Formula assures that this will be the case when both curves $C_t$ and $B$ are rational and the morphisms $B, C_t \rightarrow B_0$ are of degree two. 

 In \cite{eu} and \cite{salgado} we provide a large class of rational elliptic surfaces that contain such pencils of rational curves defined over the base field. Those are the surfaces satisfying the hypothesis of \cite[Theorem 1.3]{salgado}.

 We now state the main result of this article. It tells us that on K3 surfaces with more than one Jacobian elliptic fibration, one of the fibrations can be used to produce an independent new section in the Mordell-Weil group of the other fibration after base-change. 
 
\begin{thm}\label{t1} Let $X$ be a K3 surface defined over a number
  field $k$; suppose $X$ is endowed with at least two distinct elliptic
  fibrations $\pi_i: X \rightarrow B_i \simeq \mathbb{P}^1$, $i=1,2$. 
  
  Then there is an elliptic curve $C$ such that \[\mathrm{rank}((X \times_{B_i}C)(k(C)))
  \geq \mathrm{rank}(X(k(B_i)))+1. \]
\end{thm}

\begin{cor}\label{cor1}
 Let $X$ be a K3 surface as in Theorem \ref{t1}. If $\mathrm{rank}(X(k(B_j)))>0$, then
 \[ \# \{t\in B_i(k);r_{t}\geq r+1 \}= \infty,\]
where $j\neq i, \, r= \mathrm{rank}(X(k(B_i)))$ and $ r_t=  \mathrm{rank}(\pi_i^{-1}(t)(k))$.
\end{cor}

\section{Proof of Theorem \ref{t1}}

\subsection{Base-change}

 Let $\mathcal{E} \rightarrow B =\mathbb{P}^1$ be a Jacobian elliptic surface and $\iota: C' \hookrightarrow \mathcal{E} $ an irreducible curve. Let $\nu: C \rightarrow C'$ be its normalization. The composition of the applications above $\varphi =\nu \circ \iota \circ \pi$ gives us a finite covering $\varphi: C \rightarrow B$.
 \[
\xymatrix{ & & \mathcal{E} \ar[d]^{\pi}  \\ C \ar[r]^{\nu} \ar@/_2pc/[rr]^{\varphi} & C'
  \ar[r] \ar[ur]^{\iota} & B. } 
\]
We obtain a new elliptic surface by taking the following fibered product:
\[\pi_C: \mathcal{E}^C= \mathcal{E}\times_B C \rightarrow C \]
Each section $\sigma$ of $\mathcal{E}$ induces naturally a section in $\mathcal{E}^C$:
 $$(\sigma,id): C= B\times_B C \rightarrow \mathcal{E} \times_B C.$$ 
 
 We will call these sections $\textit{old sections}$. The surface $\mathcal{E}^C$ has also a \textit{new section} given by 
 $$\sigma^{\text{new}}_C= (\iota \circ \nu, id): C \rightarrow \mathcal{E}
\times_B C$$ 
 
 \textbf{Remark 3:}
 If $C$ is not contained in a fiber nor in a section, then the new section is different from the \textit{old} ones, but it is not necessarily linearly independent in the Mordell-Weil group.

\subsection{Main lemma}

 There are many ways of verifying if a given section is independent of a set of sections in an elliptic surface. One of them consists in computing the determinant of the height matrix associated to these sections. Since we are not interested in a specific elliptic surface we cannot compute such determinant (the height pairing requires information such as the bad fiber configuration and the generators of the Mordell-Weil group). For us the most useful way to check independence will be the following criteria:

\textbf{Criteria:} A curve $C \hookrightarrow X$ that is not a section nor a component of a fiber induces a new section in $X \times_B C$ independent of the \textit{old} ones if and only if for every section $C_0 \hookrightarrow X$ and every $n \in \mathbb{N}$, the curve $C$ is not a component of $[n]^{-1}(C_0)$.

In what follows we will need the following definition.

\begin{dfn}
Let $X$ be a surface and $\mathcal{C}$ be a set of curves in $X$. We will say that $\mathcal{C}$ is a \textit{numerical family} if all its members belong to the same numerical class on $\mathrm{Pic} X$.
\end{dfn}
 
 The following lemma states that for a numerical family of curves in $X$ all but finitely many induce an independent new section on the base-changed surface: it suffices to verify the above criteria for finitely many $n\in \mathbb{N}$ and finitely many sections $C_0$ in the Mordell-Weil group.
  
\begin{lem}\label{prpimportante} Let $\pi:X\rightarrow \mathbb{P}^1_k$ be an elliptic surface. Let 
$\mathcal{L}=\{L_t:t\in \mathbb{P}^1\}$ be a non-constant numerical family
of curves in $X$ whose
generic member is an irreducible curve which is neither a section nor
a fiber of $\pi$. Then, for almost all $t\in\mathbb{P}^1(k)$, the new section $\sigma^\mathrm{new}_{L_t}$ is independent of the old sections.  
\end{lem}

 The main tools used to prove this lemma are:

\begin{itemize}
\item[i)] Kummer theory for elliptic curves to bound the set of natural numbers $n$ that one has to consider in the criteria, and
\item[ii)] N\'{e}ron-Tate height on elliptic surfaces, which is determined by the intersection number of the curve with the components of the fibers and the generators of the Mordell-Weil group, to bound the set of sections $C_0$ such that for a fixed $n$, the curve $[n]^{-1}C_0$ contain a given curve as a component. 
\end{itemize}

 For a complete proof of it see \cite[Proposition 4.2]{salgado}.

\textbf{Proof of Theorem \ref{t1}:}
 
  By Lemma \ref{prpimportante}, it suffices to observe that the set of fibers of $\pi_j$ is a non-constant numerical family of curves in $X$ that are not contained in sections nor in the fibers of $\pi_i$ for $i \neq j$. 
  
 In fact, the curves $\pi_j^{-1}(t)$ form a numerical family since they are fibers of an elliptic fibration in $X$. Since they are (almost all) irreducible smooth genus one curves, they are not components of sections. Also, because they form a distinct fibration, by hypothesis, they cannot be components of the fibers of $\pi_i$.
 
 It follows from Lemma \ref{prpimportante} applied to the numerical family $\{\pi_j^{-1}(t); t \in B_j(k)\}$ that for all but finitely many $t \in B_j(k)$, the curves $\pi_j^{-1}(t)$ induce a new independent section on the base-changed surface. Pick $C$ to be one of these curves then $(X \times_{B_i} C)(k(C))$ has generic rank strictly larger than the generic rank of $X(k(B_i))$.
\qed

\textbf{Proof of Corollary \ref{cor1}:} 

 By the hypothesis, the new base $C= \pi_j^{-1}(t)$, for a given $t\in B_j(k)$, has infinitely many $k$-rational points. The corollary will follow from the theorem after an application of Silverman's Specialization Theorem.
\qed
 
\textbf{Remark 4:} If one of the fibrations has rank zero then the argument presented in \cite{BogoTschinkel} is enough to show that infinitely many fibers have positive rank. It suffices to find a curve $C$ such that the restriction of the fibration $C \rightarrow B$ is ramified in the smooth locus. Since the restriction of the fibration to torsion sections is \'etale the section induced by the base-change by $C$ will be of infinite order.

\section{Examples}

 \textbf{Example 1: A Kummer surface given by a double cover of a rational elliptic surface}

  We consider the K3 elliptic surface $X \rightarrow B \simeq \mathbb{P}^1$ whose
 generic fiber is given by
\[ y^2=x(x+1)(x+c^2) \,; \,\, \text{where } c=\frac{t^2-1}{t^2+1}.\]

 As remarked by Piatetski-Shapiro and Shafarevich in \cite{Sha}, this
 surface is the Jacobian of the Fermat quartic in $\mathbb{P}^3$ given by
 \[x_0^4-x_1^4=x_2^4-x_3^4 \]

 where the elliptic fibration considered is given by the equations
\[ x_0^2+x_1^2=t(x_2^2-x_3^2) \,; \, t(x_0^2-x_1^2)=x_2^2+x_3^2.\]

 Its degenerate fibers are all of type $I_4$. They are located above
 $t=, \infty, -1, 0,1,i,-i$. It follows from the Shioda-Tate formula
 applied to $X$
 that the generic rank of this elliptic fibration over $\bar{\mathbb{Q}}$ is zero.

 One sees imediately that $X$ comes quadratically from the
 rational elliptic surface $X_0$ whose generic fiber is \[
 y^2=x(x+1)(x+d^2) \,; \,\text{where } d= \frac{t-1}{t+1}\] via the
 quadratic morphism 
\[\varphi: D\simeq\mathbb{P}^1 \rightarrow B_0\simeq \mathbb{P}^1\, \, t\mapsto t^2.\]

 The singular fibers of $X_0$ are of type $I_2,I_2,I_4,I_4$ and are located above $0, \infty, 1$ and
 $-1$. The set of degenerate fibers is defined over $\mathbb{Q}$, thus the
 $I_2$ fibers, that may be conjugate under the Galois group, constribute with at least one unity to the rank of the group $\mathrm{NS}(X_0|\mathbb{Q})$, and the $I_4$ fibers contributes with at least two unities
 to the rank of $\mathrm{NS}(X_0| \mathbb{Q}$). By the Shioda-Tate formula we
 have \[ \mathrm{rank } (\mathrm{NS}(X_0| \mathbb{Q})) \geq 2+1+2=5.\]

  Thus a $k$-minimal model of $X_0$ has degree at least three. It follows from \cite[Corollary 0.2.4]{salgado} that $X_0$ has infinitely many fibers with rank equal to at least the generic rank plus one. Since $X$ dominates $X_0$ the same holds true for $X$. 

\textbf{Example 2: Two fibrations on a Kummer surface not coming quadratically from a rational elliptic surface:}
  
 In \cite{Nishiyama} Nishiyama describes all Jacobian fibrations associated to certain Kummer surfaces. We consider here one of these surfaces. 

 Let $E_{\zeta_3}= \mathbb{C}/(\mathbb{Z} +\zeta_3 \mathbb{Z})$ and $X=X_3$,
 defined over $k=\mathbb{Q}(\zeta_3)$, the minimal resolution of the surface \[ E_{\zeta_3}
 \times E_{\zeta_3}/ <\sigma: (z_1,z_2)\mapsto (\zeta_3 z_1,
 \zeta_3^{-1} z_2)> ,\, \, (\zeta_3= \frac{-1+\sqrt{-3}}{2}).\]

 The surface $X$ is K3 and is endowed with six distinct Jacobian fibrations \cite[Theorem 3.1]{Nishiyama} among which two have positive generic rank over $k$.

 The two fibrations are the following:

\[\pi_1:X \rightarrow B_1 ;\, \text{ with bad fibers of type } 2III^*, I_6^*,\]
 and
\[ \pi_2:X \rightarrow B_2 ; \, \text{ with bad fibers of type } III^*, I_{18}.\]

 The fibration $\pi_1$ cannot come quadratically from a rational elliptic surface since such a surface can never have fibers of type $I_6^*$, and a fiber of type $I_3^*$ induces either two fibers of
 type $I_3^*$, or one fiber de type $I_6$ after a quadratic base-change. The same holds for the fibration $\pi_2$ since no fiber becomes  a fiber of type $III^*$ after a quadratic base-change.

 The curves in the pencil \[ \mathcal{L}=\{ C_t = \pi_2^{-1}(t);
 t \in \mathbb{P}^1(k) \} \]

 are not components of the fibers nor of the sections of $\pi_1$. They satisfy: \[\# C_t(k)=
 \infty \text{ for almost all } t \in \mathbb{P}^1(k).\] It follows from Lemma \ref{prpimportante} that almost all $C_t \in \mathcal{L}$ induce an independent section in $X \times_{B_1} C_t$. We have

\[\mathrm{rank}(X \times_{B_1} C_t(k(C_t)))\geq
\mathrm{rank}(X(k(B_1)))+1 \]
and
\[\# \{t \in B_1(k); r_t \geq r_1+1 \}= \infty, \]
where $r_1=1$ is the generic rank of $\pi_1$ and $r_t$ is the rank of the Mordell-Weil group $\pi_1^{-1}(t)(k)$.

\section*{Acknowledgements}{ The author would like to thank Marc Hindry for several discussions on the topic, and the referee for comments and suggestions that improved the readability of the paper.}

\end{document}